\documentclass{ifacconf}
\usepackage{amsmath}
\usepackage{amsfonts}
\usepackage{natbib} 
\usepackage{amssymb}
\usepackage{graphicx}
\usepackage{color}
\usepackage{dsfont}
\usepackage{mathrsfs}
\usepackage[usenames,dvipsnames]{pstricks}
\usepackage{enumerate} 
\usepackage{easyReview} 

\allowdisplaybreaks
\usepackage{flushend}
\newcommand{\rev}{\textcolor{black}}
\newcommand{\modif}{\textcolor{black}}

\begin{document}
\begin{frontmatter}


\title{Output-feedback stabilization of a class of $n+m$ linear hyperbolic ODE-PDE-ODE systems}

\author[First]{Jean Auriol,} \author[Second]{Federico Bribiesca Argomedo,} 
\address[First]{Universit\'e Paris-Saclay, CNRS, CentraleSup\'elec, Laboratoire des Signaux et Syst\`emes, 91190, Gif-sur-Yvette, France. }
\address[Second]{INSA Lyon, Universite Claude Bernard Lyon 1, Ecole Centrale de Lyon, CNRS, Amp\`ere, UMR5005, 69621 Villeurbanne, France}

\begin{abstract}
In this paper, we design an output-feedback controller to stabilize $n + m$ hetero-directional transport partial differential equations (PDEs) coupled on both domain boundaries to ordinary differential equations (ODEs). This class of systems can represent, for instance, actuator and load dynamics at the boundaries of a hyperbolic system. The actuator is located at the connection point between the PDE and one of the ODEs, and we consider anti-collocated PDE measurements. We first design a state-observer by combining the backstepping methodology with time-delay system approaches. We then introduce a state feedback controller using analogous techniques before designing the wanted output-feedback control law.

\end{abstract}
\end{frontmatter}
\vspace{-0.1cm}

\section{Introduction}

Several recent contributions have focused on the control of interconnected systems encompassing hyperbolic partial differential equations (PDEs) and ordinary differential equations (ODEs), given their prevalence in diverse industrial processes, such as electric power transmission systems~\citep{schmuck2014feed} and traffic networks~\citep{yu2023traffic}. This category of networks is particularly adept at modeling complex industrial processes as the propagation of torsional waves in drilling systems~\citep{aarsnes2018torsional} or deepwater construction vessels~\citep{stensgaard2010subsea}. These interconnected networks are commonly referred to as ODE-PDE-ODE systems, where the ODEs capture the dynamics of actuators and loads.

Most constructive control methodologies for interconnected PDE-ODE systems predominantly adopt the backstepping approach.  In the seminal paper~\citep{krstic2008backstepping}, a re-interpretation of the classical Finite Spectrum Assignment \citep{Manitius1979} was proposed, modeling ODEs with input delays as PDE-ODE interconnections. The consequential impact of this innovative approach has been far-reaching, enabling the formulation of methodologies for designing observers, controllers, and parameter estimation methods across a diverse spectrum of interconnected systems. Noteworthy applications include systems with varying delays~\citep{bekiaris2013nonlinear,bresch2012commande} and cascades of PDEs~\citep{auriol2019interconnected,redaud2021SCL}, or cascaded interconnections of hyperbolic PDE-ODE systems, such as~\citep{aamo2012disturbance,deutscher2021backstepping,auriol2018delay,hasan2016boundary,zhou2012}. More recently, constructive results have emerged for non-cascaded PDE-ODE systems. In particular, a stabilizing state-feedback control law has been proposed in~\citep{di2018stabilization,wang2018control}. A PDE-ODE-PDE configuration was considered in~\citep{auriol2020robust_bis}.

For ODE-PDE-ODE configurations, an output-feedback controller has been designed in~\citep{deutscher2018output} based on assumptions that guarantee the existence of a Byrnes–Isidori normal form for one of the ODE. These restrictions are partially avoided in~\citep{BouSaba2019}, where the control design relies on a \emph{rewriting of the interconnection as a time-delay system}. This approach was later extended in \citep{wang2020delay,auriol2022observer,redaud2024output} to encompass a state observer. Some recent developments have also been obtained for interconnected PDE systems with non-linear ODEs using a modular design of tracking controllers~\citep{irscheid2021observer}. In all these contributions, the authors assumed the control input was acting on the ODE state. The cases where the control input acts at the connection point between the PDE and one of the ODEs have not been well studied in the literature, and most of the contributions neglect the actuator dynamics in such a configuration~\citep{auriol2018delay,auriol2020robust,de2018backstepping}. In~\citep{deutscher2018output}, a state observer was designed for an ODE-PDE-ODE system where the measurement corresponds to the PDE state. In the case of a $2\times 2$ PDE system, a stabilizing controller was designed in~\citep{auriol2023robustification}.

In this paper, we extend these results to design a stabilizing output-feedback controller for a system of $n+m$ linear first-order hyperbolic Partial Differential Equations coupled with Ordinary Differential Equations at both boundaries of a one-dimensional spatial domain. The control input acts at one of the PDE boundaries, and the available PDE measurements are anti-collocated. The presented approach expands upon the methodology introduced in~\citep{auriol2023robustification} by combining the backstepping technique with time-delay approaches. Using an invertible integral transformation, we map the system into a simpler target system. Then, we design a state observer for this target system using a time-delay representation. A state-feedback controller is obtained using analogous techniques. Finally, the output-feedback controller is obtained by combining the state-feedback controller with the previously designed observer after adding an adequate low-pass filter to guarantee the existence of robustness margins~\citep{auriol2023robustification}. 

\emph{Notations}
The state space is $\chi = \mathbb{R}^{p} \times L^2([0,1]; \mathbb{R})^{n+m} \times \mathbb R^{q}$, where $p, n, m, q$ are positive integers.
For \modif{$(X_0,u,v,X_1) \in \chi$}, we introduce the corresponding $\chi-$norm $$ \modif{||(X_0,u,v,X_1)||^2_\chi=||X_0||^2_{\mathbb{R}^p}+||u||^2_{L^2}+||v||^2_{L^2}+||X_1||^2_{\mathbb{R}^q}}.$$ We denote $s$, the Laplace variable.
\section{Problem statement} \label{Sec_prob}
\subsection{Presentation of the system}
In this paper, we consider a \modif{$n + m$} linear hetero-directional hyperbolic system coupled through its boundaries with linear ODEs: 
\begin{eqnarray}
\dot{X}_0(t)&=& A_0X_0(t)+E_0v(t,0), \label{eq_X_0}\\
u(t,0)&=&C_0X_0(t)+Q v(t,0)+U(t) , \\
u_t+\Lambda^+ u_x &=& \Sigma^{++}(x) u(t,x)+ \Sigma^{+-}(x) v(t,x), \label{eq_u}\\
v_t-\Lambda^- v_x &=& \Sigma^{-+}(x) u(t,x)+ \Sigma^{--}(x) v(t,x), \label{eq_v}\\
v(t,1)&=&Ru(t,1)+C_1X_1(t), \label{eq_bound} \\
\dot{X}_1(t)&=& A_1X_1(t)+E_1u(t,1), \label{eq_ODE_distal}
\end{eqnarray}
defined for a.e. $(t,x) \in [0, +\infty) \times [0,1]$. The state of the system is $(X_0(t), u(t, \cdot), v(t, \cdot), X_1(t)) \in \chi$.
 The initial condition is taken as $((X_0)_0, u_0, v_0, (X_1)_0) \in \chi$ and we consider weak solutions to \eqref{eq_X_0}-\eqref{eq_ODE_distal} \citep{bastin2016stability}. The open-loop system is well-posed in the sense of ~\citep[Theorem A.6, page 254]{bastin2016stability}. The matrices $\Lambda^+$ and $\Lambda^-$ are diagonal and represent the transport velocities. We have $\Lambda^+=$ diag $(\lambda_i)$ and $\Lambda^-=$ diag ($\mu_i$) and we assume that their coefficients satisfy
$
-\mu_m < \cdots< -\mu_1 <0<\lambda_1 < \cdots < \lambda_n.
$
{The case of equal transport velocities or zero velocities can be overcome under additional assumptions following the methodology presented in~\citep{de2024backstepping,chen2023block}.}
The spatially-varying matrices $\Sigma^{\cdot \cdot}$ are \modif{continuous} (each coefficient of the matrix is a continuous function). \modif{With no loss of generality, we assume that the matrices $\Sigma^{++}$ and $\Sigma^{--}$ have zero diagonal elements \rev{ \citep{hu2015boundary}}}. The different coupling matrices satisfy $A_0 \in \mathbb{R}^{p \times p}$, $E_0 \in \mathbb{R}^{p \times m}$, $C_0 \in \mathbb{R}^{n \times p}$, $A_1 \in \mathbb{R}^{q \times q}$, $E_1 \in \mathbb{R}^{q \times n}$, $C_1 \in \mathbb{R}^{m \times q}$, $R \in \mathbb{R}^{m \times n}$, $Q \in \mathbb{R}^{n \times m}$. The control input $U(t)$ belongs to $\mathbb{R}^{n}$.
We consider the case of anti-collocated measurement, i.e., the measurement $y(t)$ is defined by
\begin{align}
    y(t)\doteq u(t,1). \label{measurement}
\end{align}
Due to the \rev{symmetry} of the system, measuring the $v(t,0)$-state would not change the nature of the problem. Compared to~\citep{deutscher2018output}, we consider that the control input is located at the junction between the ODE $X_0$ and the PDE. This configuration has not been well-studied in the literature and requires specific control approaches.

\begin{figure}[htb]%
\begin{center}
	\begin{tikzpicture}
		\draw [>=stealth,color=green!50!black!90,very thick] (-1,0) -- (0,0)--(0,-1.5)--(-1,-1.5)--(-1,0);
		\draw [color=green!50!black!90] (-0.5,-0.7) node[above]{$X_0$};
  \draw [color=green!50!black!90] (-0.5,-1.3) node[above]{ODE};
\draw [>=stealth,->,color=green!50!black!90,very thick] (0.1,0) -- (0.5,0);
\draw [color=green!50!black!90] (0.2,0) node[above]{$C_0$};
\draw [>=stealth,<-,color=green!50!black!90,very thick] (0.1,-1.5) -- (0.5,-1.5);
\draw [color=green!50!black!90] (0.2,-1.5) node[below]{$E_0$};

\draw [>=stealth,<-,gray,very thick] (2.5,0.3) -- (2.5,0);
\draw [gray] (2.5,0.3) node[above]{$y(t)$};
\draw [>=stealth,->,purple,very thick] (0.75,0.3) -- (0.75,0);
\draw [purple] (0.75,0.3) node[above]{$U(t)$};
\draw [>=stealth,->,red,very thick] (0.75,0) -- (2.5,0);
\draw [red] (1.6,0) node[above]{$u(t,x)$};
\draw [>=stealth,<-,blue,very thick] (0.75,-1.5) -- (2.5,-1.5);
\draw [blue] (1.6,-1.5) node[below]{$v(t,x)$};
\draw [>=stealth,<-,dashed, thick] (1.2,-1.5) -- (1.2,0);
\draw [>=stealth,->,dashed, thick] (1.8,-1.5) -- (1.8,0);
\draw [red,>=stealth, thick](2.9,-0.75) arc (0:45:1.1);
\draw [red,>=stealth,->, thick](2.9,-0.75) arc (0:-45:1.1);
\draw [red] (2.9,-0.75) node[left]{$R$};
\draw [blue,>=stealth, thick](0.35,-0.75) arc (-180:-135:1.1);
\draw [blue,>=stealth,->, thick](0.35,-0.75) arc (-180:-225:1.1);
\draw [blue] (0.35,-0.75) node[right]{$Q$};
\draw [>=stealth,<-,orange,very thick] (2.8,-1.5) -- (3.2,-1.5);
\draw [orange] (3,-1.5) node[below]{$C_1$};
\draw [>=stealth,->,orange,very thick] (2.8,0) -- (3.2,0);
\draw [orange] (3,0) node[above]{$E_1$};
\draw [>=stealth,color=orange,very thick] (3.3,0) -- (4.3,0)--(4.3,-1.5)--(3.3,-1.5)--(3.3,0);
		\draw [color=orange] (3.8,-0.7) node[above]{$\dot{X}_1$};
  \draw [color=orange] (3.8,-1.3) node[above]{ODE};
\draw [>=stealth,->,black,very thick] (0.6,-2.3) -- (2.9,-2.3);
\draw [>=stealth,black,very thick] (0.9,-2.4) -- (0.9,-2.2);
\draw [>=stealth,black,very thick] (2.5,-2.4) -- (2.5,-2.2);
\draw [black] (0.9,-2.4) node[below]{0};
\draw [black] (2.5,-2.4) node[below]{1};
\draw [black] (2.9,-2.3) node[right]{$x$};
\end{tikzpicture}
\end{center}
\caption{Schematic representation of the system~\eqref{eq_X_0}-\eqref{eq_ODE_distal}. }%
\label{Fig_Example_ODE_PDE_ODE}%
\end{figure}
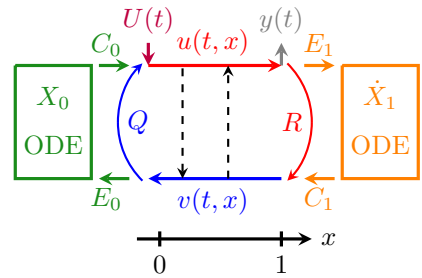

\subsection{General objectives and assumptions}
The objective of this paper is to design an output feedback controller for the system~\eqref{eq_X_0}-\eqref{eq_ODE_distal} based on the available measurement $y(t)$ given by equation~\eqref{measurement}. 
\modif{Even if the velocities are considered constant here, the proposed methodology could be extended to the case of spatially-varying transport velocities (following the approach given in \citep{hu2015boundary}).} 
We now make the following assumption
\begin{assum} \label{Assum_OL_stab}
The system defined for all $i \in [1,n]$ by 
\begin{align}
   z(t)=\sum_{k=1}^m\sum_{\ell=1}^n Q_{ik}R_{k\ell}  z(t-\frac{1}{\mu_k}-\frac{1}{\lambda_\ell}), \label{eq_assum_OL}
\end{align}
is exponentially stable.
\end{assum}
 Assumption~\ref{Assum_OL_stab} constitutes a reasonable assumption since it prevents system~\eqref{eq_X_0}-\eqref{eq_ODE_distal} from having an asymptotic chain of eigenvalues with non-negative real parts \citep{halebook,auriol2019explicit}. It has been shown in \citep{logemann1996conditions} that having an open-loop transfer function with an infinite number of poles on the closed right half-plane implies no (delay-)robustness margins in closed-loop (i.e., the introduction of any arbitrarily small delay in the actuation will destabilize the closed-loop system). Therefore, Assumption~\ref{Assum_OL_stab} is slightly stronger than a necessary condition for delay-robust stabilization. If the delays are rationally independent\footnote{\rev{Extending the variable $z$, it is always possible to rewrite the system in a form with either a single delay, or where all delays are rationally independent~\citep[Chapter 9]{halebook}.}}, Assumption~\ref{Assum_OL_stab} is equivalent to the following condition~\citep{halebook}: $
	\underset{\theta_{k\ell} \in [0,2\pi]^{n\times m}}{\sup} \text{Sp}~ (\sum_{k=1}^m\sum_{\ell=1}^n Q_{ik}R_{k\ell}  \exp(\modif{j} \theta_{k\ell}))<1,$
where Sp denotes the spectral radius, \modif{$j$ is the imaginary unit}, \rev{and the $\theta_{k\ell}$ belong to  $[0,2\pi]^{n\times m}$}. This condition is simplified if the delays are rationally dependent. Furthermore, since the spectral radius of a matrix is upper-bounded by any norm of the matrix, easy to compute sufficient conditions for this spectral radius condition to hold can be derived using different norms of the matrices involved at the cost of increased conservatism. Similarly to what has been done in~\citep{auriol2023robustification}, we will need additional stabilizability and detectability conditions (Assumptions~\ref{Chap_2_ass_detec_PDE} and~\ref{Chap_2_ass_stab_PDE}) that will be given later. 

\section{Observer design} \label{Sec_observer}

This section is dedicated to designing a state observer for the system governed by equations \eqref{eq_X_0}-\eqref{eq_ODE_distal}, leveraging the available measurement given in equation \eqref{measurement}.
To simplify the computations and the design of the observer, we will work with a target system obtained from \eqref{eq_X_0}-\eqref{eq_ODE_distal} using a backstepping transformation. This transformation relocates the in-domain coupling terms, denoted as $\Sigma^{\cdot \cdot}$, to the actuated boundary.

\subsection{Backstepping transformation} \label{Sec_transfo}
Inspired by~\citep{hu2015boundary}, we use an integral transformation to move the local coupling terms~$\Sigma^{\cdot \cdot}$ to the actuated boundary (in the form of integral terms). 
Consider the Volterra  transformation $\mathcal{T}$, similar to the one introduced in~\citep{hu2015boundary,auriol2019explicit,auriol2022observer}
\begin{eqnarray}
X_0(t)&=&\xi(t)-\int_0^1 L_1(\nu) \alpha(\nu)+L_2(\nu) \beta(\nu) dy,\\ \label{back_X0_ODE}
u(t,x)&=&\alpha(t,x)-\int_x^1 L^{\alpha \alpha}(x,\nu) \alpha(\nu)d\nu \nonumber \\
&-&\int_x^1 L^{\alpha \beta}(x,\nu) \beta(\nu)d\nu +\gamma_\alpha(x)X_1(t), \label{back_1} \\
v(t,x)&=&\beta(t,x)-\int_x^1 L^{\beta \alpha}(x,\nu) \alpha(\nu)dy, \nonumber \\
&-&\int_x^1 L^{\beta \beta}(x,\nu) \beta(\nu)d\nu +\gamma_\beta(x)X_1(t), \\
X_1(t)&=&X_1(t),\label{back_2}
\end{eqnarray}
where the kernels are bounded functions defined either on $\mathcal{T}_u=\{(x,z) \in [0,1]^2,~x\leq z \}$, or $[0,1]$.
This transformation rewrites $(X_0, u,v, X_1)=\mathcal{T}(\xi, \alpha, \beta,X_1)$. \modif{Denoting $\Lambda = \text{diag} (\Lambda^+, -\Lambda-)$, $\Sigma = \left(\begin{smallmatrix} \Sigma^{++}&\Sigma^{+-}\\\Sigma^{-+}&\Sigma^{--}\end{smallmatrix}\right)$ and $L = \left(\begin{smallmatrix} L^{\alpha \alpha}&L^{\alpha \beta}\\L^{\beta \alpha}&L^{\beta \beta}\end{smallmatrix}\right)$, $\gamma=(\gamma_\alpha, \gamma_\beta)$, we obtain}
\modif{
\begin{align}
    &\Lambda L_x + L_\nu \Lambda = \Sigma(x) L \label{eq_kernel_L_aa},\quad
    \Lambda \gamma_x(x) =\Sigma(x) \gamma - \gamma^\top  A_1, \\
    &(L_1(x))_x\Lambda^+= A_0 L_1(x)+E_0L^{\beta \alpha}(0,x), \label{eq_kernel_L1}\\
    &(L_2(x))_x\Lambda^-=-A_0 L_2(x)-E_0L^{\beta \beta}(0,x),\label{eq_kernel_L2}
\end{align}
with the boundary conditions
\begin{align}
    &\Lambda L(x,x) - L(x,x)\Lambda = \Sigma(x),\\
    &  L_1(0)=0,~L_2(0)\Lambda^-=L_1(0)\Lambda^+Q+E_0,\label{eq_bound_L}
\end{align} and $\gamma_\alpha(1)=0, ~ \gamma_\beta(1)=C_1$.} Finally, we define $L_{ij}^{\alpha \alpha}(0,\nu)$ for $i \leq j$ by
\begin{align}
    L^{\alpha \alpha}_{ij}(0,\nu)=(QL^{\beta \alpha}(0,\nu)+C_0L_1(\nu))_{ij}.\label{eq_kernel_final}
\end{align}
 To this set of equations, we add arbitrary values for $L_{ij}^{\alpha \alpha}(x,1)$ (when $i>j$) and $L_{ij}^{\beta \beta}(x,1)$ (when $i>j$) and $L_{ij}^{\beta \beta}(0,\nu)$ (when $i\leq j$).
\modif{Reinterpreting the ODEs in \eqref{eq_kernel_L1}-\eqref{eq_kernel_L2} as PDEs evolving in the triangular domain $\mathcal{T}_u$ with horizontal characteristic lines (since there is only
an evolution along the $x$ axis), it is possible to adjust the results from \citep[Theorem 3.2]{di2018stabilization} to guarantee that the set of PDEs and ODEs \eqref{eq_kernel_L_aa}-\eqref{eq_kernel_final} has a unique solution which is piecewise continuous.} 
 The boundedness of transformation \eqref{back_X0_ODE}-\eqref{back_2} is a direct consequence of the structure of the transform (identities, integral operator and matrices) and the regularity of the different kernels. Its invertibility is a consequence of the structure of the transformation, which is block triangular with the blocks on the diagonal being either identities (for the ODEs) or invertible Volterra operators (for the PDEs).
 \subsection{Target system}
The invertible backstepping transformation~\eqref{back_X0_ODE}-\eqref{back_2} maps the original system~\eqref{eq_X_0}-\eqref{eq_ODE_distal}  to the following target system 
\begin{align}
&\dot{\xi}(t)= A_0 \xi(t)+G_3\alpha(t,1)+G_4X_1(t), \label{eq_xi_target}\\
&\alpha(t,0)=Q\beta(t,0)+C_0 \xi(t)+(Q\gamma_\beta(0)-\gamma_\alpha(0))X_1(t) \nonumber \\
&+\int_0^1 F^\alpha(\nu)\alpha(t,\nu)+F^\beta(\nu) \beta(t,\nu)d\nu+U(t),\\
 &\alpha_t(t,x)+\Lambda^+ \alpha_x(t,x) = G_1(x)\alpha(t,1) , \label{eq_alpha_obs}\\
&\beta_t(t,x)-\Lambda^- \beta_x(t,x) = G_2(x)\alpha(t,1), \label{eq_beta_obs}\\
&\beta(t,1)=R\alpha(t,1),~\dot{X}_1(t)=  A_1X_1(t)+E_1\alpha(t,1). \label{eq_ODE_distal_target_obs}
\end{align} 
The functions $G_1$ and $G_2$ satisfy
\begin{align}
    &G_1(x)=\int_x^1 L^{\alpha\alpha}(x,\nu)G_1(\nu)+L^{\alpha\beta}(x,\nu)G_2(\nu)d\nu \nonumber \\
    &-L^{\alpha\alpha}(x,1)\Lambda^++L^{\alpha \beta}(x,1)\Lambda^-R-\gamma_\alpha(x)E_1, \label{eq_G_0} \\
    &G_2(x)=\int_x^1 L^{\beta\alpha}(x,\nu)G_1(\nu)+L^{\beta\beta}(x,\nu)G_2(\nu)d\nu \nonumber \\
    &-L^{\beta\alpha}(x,1)\Lambda^++L^{\beta \beta}(x,1)\Lambda^-R-\gamma_\beta(x)E_1. \label{eq_G_1}
\end{align}
The set of equations \eqref{eq_G_0}-\eqref{eq_G_1} has a unique solution (Volterra equations of the second kind~\citep{yoshida1960lectures}). 
The matrices $G_3$ and $G_4$ are defined by
    $G_3(x)=L_2(x,1)\Lambda^-R-L_1(x,1)\Lambda^+
    +\int_0^1(L_1(x)G_1(x)+L_2(x)G_2(x))dx,$ and $
    G_4(x)=E_0\gamma_\beta(0).$
Finally, the matrix $F^\beta$ the matrix $F^\alpha$ are defined by
$F^\alpha(\nu)=L^{\alpha \alpha}(0,\nu)-QL^{\beta \alpha}(0,\nu)-C_0L_1(\nu),$ and $
    F^\beta(\nu)=L^{\alpha \beta}(0,\nu)-QL^{\beta \beta}(0,\nu)-C_0L_2(\nu).$
Note that $F^\alpha$ is \textbf{strictly lower triangular} due to equation \eqref{eq_kernel_final}. The measurement $y(t)$ remains unchanged.
\subsection{Time-delay formulation}
 Applying the method of characteristics (see~\citep{auriol2019explicit,auriol2023robustification}), we can show that for all $t>\tau=\frac{1}{\lambda_1}+\frac{1}{\mu_1}$
\begin{align}
    &\alpha(t,1)=\sum_{i=1}^N F_i^\alpha  \alpha(t-\tau_i,1)+\sum_{i=1}^n F_i^\xi \xi(t-\frac{1}{\lambda}_i) \nonumber \\
    &+\sum_{i=1}^n F_i^X  X_1(t-\frac{1}{\lambda}_i)+\int_0^\tau  H(\nu) \alpha(t-\nu,1)d\nu \nonumber \\
    &+\sum_{i=1}^n F_i^UU(t-\frac{1}{\lambda}_i), \label{Chap_2_eq_delay_alpha_1}
\end{align}
where $N\in\mathbb{N}$, and where the $\tau_i\leq \tau$ are positive delays that depend on the $\lambda_i$ and $\mu_i$, the matrices~$ F^\alpha_i$,~$F^\xi_i$,~$ F^X_i$, $F^X_i, F_i^U$, and~$ H$ can be explicitly computed following the methodology of \citep{auriol2019explicit}. 
Since the functions $U(t)$ and $y(t)=\alpha(t,1)$ are known, we can consider $y_1$ as an available measurement, where $y_1$ is defined by \begin{align}y_1(t)&=y(t)-\sum_{i=1}^NF_i^\alpha y(t-\tau_i)-\int_0^\tau H(\nu) y(t-\nu)d\nu \nonumber \\
&-\sum_{i=1}^nF_i^UU(t-\frac{1}{\lambda_i})\nonumber \\
&=\sum_{i=1}^n F_i^X \xi(t-\frac{1}{\lambda}_i)+ F_i^X  X_1(t-\frac{1}{\lambda}_i).
\end{align}
Define $Z$ as the concatenation of the vectors $\xi$ and $X_1$ and introduce the matrix $A_o=\begin{pmatrix}
    A_0 & G_4 \\ 0_{q\times p} & A_1
\end{pmatrix}$. We obtain
\begin{align}
    \dot{Z}(t)=A_0Z(t)+G_Z\alpha(t,1), \label{Chap_2eq_Z}
\end{align}
where $G_Z=\begin{pmatrix}G_3 \\E_1
\end{pmatrix}$. Defining $F_i^o=\begin{pmatrix}
    F_i^\xi & F_i^X
\end{pmatrix}$, we have $y_1(t)=\sum_{i=1}^n F_i^oZ(t-\frac{1}{\lambda_i})$. Applying Duhamel's formula to equation~\eqref{Chap_2eq_Z}, we obtain for $t>\frac{1}{\lambda_1}$
\begin{align*}
    y_1(t)=\sum_{i=1}^nF_i^o  \mathrm{e}^{-\frac{A_o}{\lambda_i}}\big(Z(t)-\int_{t-\frac{1}{\lambda_i}}^t\mathrm{e}^{A_o(t-\nu)}G_Zy(\nu)d\nu \big).
\end{align*}
We make the following assumption 
\begin{assum} \label{Chap_2_ass_detec_PDE}
  There exists~$L_o$ such that the matrix~$A_o+L_o\sum_{i=1}^nF^o_i\mathrm{e}^{-\frac{A_o}{\lambda_i}}$ is Hurwitz.
\end{assum}
This assumption can be seen as a detectability condition necessary to reconstruct the ODE states. 

\subsection{Observer equations}

We can now design the wanted state-observer for the target system \eqref{eq_xi_target}-\eqref{eq_ODE_distal_target_obs}. Herein, the notation $\hat \cdot$ is introduced as a superscript to denote the estimated states, while $\tilde \cdot$ represents the error state, defined as the difference between the actual state and the observer state. The primary objective is to achieve the convergence of the estimated state to the real state, or equivalently, the convergence of the error state to zero, with respect to the $\chi$-norm. Subsequently, a state observer is constructed to replicate these dynamics, incorporating specific output injection terms. The observer state $(\hat \xi, \hat \alpha, \hat \beta, \hat X_1)$ (or equivalently $(\hat Z, \hat \alpha, \hat \beta)$, where $\hat Z$ as the concatenation of $\hat {\xi}$ and $\hat X_1$) is the solution of a set of equations that is a copy of the original dynamics to which we add dynamical output injection gains.  
We denote $\tilde y(t)=y(t)-\hat \alpha(t,1)$, the difference between the real output and the observer output. 
The observer equations read as 
\begin{align}
&\dot{\hat{Z}}(t)= A_o \hat Z(t)+G_Z y(t)-L_o(y_1(t)-\sum_{i=1}^nF_i^o\mathrm{e}^{-\frac{A_o}{\lambda_i}}\hat Z(t))\nonumber \\
&\quad \quad \quad-L_o\sum_{i=1}^nF_i^o\mathrm{e}^{-\frac{A_o}{\lambda_i}}\int_{t-\frac{1}{\lambda_i}}^t\mathrm{e}^{A_o(t-\nu)}G_Zy(\nu)d\nu,\label{Chap_2_eq_xi_target_estim_bis}\\
&\hat\alpha(t,0)=Q\hat\beta(t,0)+C_0 \hat{\xi}(t)+(Q\gamma_\beta(0)-\gamma_\alpha(0))\hat X_1(t)\nonumber\\
&+\int_0^1 F^\alpha(\nu)\hat\alpha(t,\nu)+F^\beta(\nu)\hat \beta(t,\nu) dy-\mathcal{O}_0(\tilde y(t))\nonumber \\
&+U(t),\\
&\partial_t\hat\alpha(t,x)+\Lambda^+\partial_x\hat\alpha(t,x) = G_1(x)y(t), \\
&\partial_t\hat\beta(t,x)-\Lambda^- \partial_x\hat\beta(t,x) = G_2(x)y(t), \\
&\hat \beta(t,1)=R\hat \alpha(t,1), 
\label{Chap_2_eq_xi_target_estim_bis_end}
\end{align}
with any (arbitrary) initial conditions in $\chi$.
 The operator $\mathcal{O}_o$ still has to be defined.
The error system is obtained by subtracting the observer dynamics from the real one. We obtain
\begin{align}
&\dot{\tilde Z}(t)= A_o \tilde Z(t)+L_o\sum_{i=1}^nF^o_i\mathrm{e}^{-\frac{A_o}{\lambda_i}}\tilde Z(t), \\
&\tilde \alpha(t,0)=C_0\tilde \xi(t)+Q\tilde \beta(t,0)+(Q\gamma_\beta(0)-\gamma_\alpha(0))\tilde X_1(t) \nonumber \\
&+\int_0^1 F^\alpha(\nu)\tilde\alpha(t,\nu)+F^\beta(\nu)\tilde\beta(t,\nu)d\nu +\mathcal O_0( \tilde y(t)), \\
&\partial_t\tilde\alpha(t,x)+\Lambda^+\partial_x\tilde  \alpha(t,x) = 0, \label{eq_alpha_error}\\
&\partial_t\tilde \beta(t,x)-\Lambda^- \partial_x\tilde \beta(t,x) =0, \label{eq_beta_error} \\
&\tilde\beta(t,1)=R\tilde\alpha(t,1).
 \label{Chap_2_eq_X_1_error_bis}
\end{align}
 To guarantee the exponential stability of the error system, it is sufficient to show the convergence of $\tilde \xi$, $\tilde \alpha(t,1)$ and $\tilde X_1$ to zero. More precisely, we have the following lemma
\begin{lem}\label{lem_conv_cascade}
If $\tilde \xi(t)$, $\tilde \alpha(t,1)$ and $\tilde X_1(t)$ exponentially converge to zero, then the state $(\tilde \xi, \tilde \alpha, \tilde \beta, \tilde X_1)$ converges to zero in the sense of the $\chi$-norm. This implies the convergence of the observer state to the real state.
\end{lem}
\begin{pf}
Due to the stability of the observer operators and using the transport structure of \eqref{eq_alpha_error} and \eqref{eq_beta_error}, the exponential convergence of $\tilde X_1$ and $\tilde \alpha(t,1)$ to zero imply the exponential convergence of the states $\tilde \alpha(t,x)$ and $\tilde \beta(t,x)$.  \end{pf}
\subsection{Design of the operator $\mathcal{O}_0$}
We now want to define the operator $\mathcal{O}_0$ such that $\tilde \xi$, $\tilde \alpha(t,1)$ and $\tilde X_1$ exponentially converge to zero.  Since $A_o+L_o\sum_{i=1}^nF^o_i\mathrm{e}^{-\frac{A_o}{\lambda_i}}$ is Hurwitz due to Assumption~\ref{Chap_2_ass_detec_PDE}, we already have the exponential stability of $\tilde \xi$ and $\tilde X_1$. Following the methodology proposed in \citep{auriol2019explicit}, we can apply the method of characteristics to obtain for $t>\tau$
\begin{align}
    &\tilde \alpha(t,1)=\sum_{i=1}^N F_i^\alpha \tilde \alpha(t-\tau_i,1)+\sum_{i=1}^n F_i^\xi \tilde \xi(t-\frac{1}{\lambda}_i)+\bar F_i^X \nonumber \\
    &\tilde X_1(t-\frac{1}{\lambda}_i) +\int_0^\tau \bar H(\nu)\tilde \alpha(t-\nu,1)d\nu+\mathcal{O}_0(\tilde y(t)),
\end{align}
where the matrices~$F^\alpha_i$,~$F^\xi_i$,~$\bar F^X_i$ and~$\bar H$ are identical to the ones given in equation~\eqref{Chap_2_eq_delay_alpha_1}. Thus, we choose $\mathcal{O}_0(\tilde y(t))$ as
\begin{align}
   \mathcal{O}_0(\tilde y(t))=&-\sum_{i=1}^N F_i^\alpha\tilde y(t-\tau_i) \nonumber \\
   &-\int_0^\tau H(\nu)\tilde y(t-\nu)d\nu.\label{Chap_2_operator_O_0}
\end{align}
We can now write the following theorem
\begin{thm} \label{Chap_2_th_observer_PDE}
Consider the operators  $\mathcal{O}_0$ defined by equation~\eqref{Chap_2_operator_O_0}. Consider that Assumption~\ref{Assum_OL_stab} and Assumption~\ref{Chap_2_ass_detec_PDE} are verified. Define the observer state $(\hat{X}_0, \hat u,  \hat v, \hat{X}_1)=\mathcal{T}( \hat{\xi},$ $ \hat \alpha, \hat \beta,\hat X_1)$, where $(\hat{\xi}, \hat \alpha, \hat \beta,\hat X_1)$ is the solution of the system~\eqref{Chap_2_eq_xi_target_estim_bis}-\eqref{Chap_2_eq_xi_target_estim_bis_end}. Then the state $(\hat X_0, \hat u,  \hat v, \hat X_1)$ exponentially converges to $(X_0, u,  v, X_1)$ in the sense of the $\chi$-norm.
\end{thm}
\begin{pf}
 We have already shown that $\tilde X_1$ and $\tilde \xi$ exponentially converge to zero.  With this choice of operator, $\tilde \alpha(t,1)$ exponentially converges to zero. Consequently, Lemma~\ref{lem_conv_cascade} implies that the state $(\tilde \xi_0, \tilde \alpha, \tilde \beta, \tilde X_1)$ exponentially converges to zero for the $\chi$-norm. Using the invertibility and boundedness of the linear transformation $\mathcal{T}$, we conclude the proof.
\end{pf}
 Note that the observer operator $\mathcal{O}_0$ may not be strictly proper \rev{(due to the cancellation of the boundary reflection terms)} and the observer system may consequently be sensitive to delays in the measurements. \rev{Low-pass filtering of the measurement will lead to a strictly proper observer operator while still guaranteeing the convergence of the estimated states towards the real states. This will consequently allow robustness margins to exist~\citep{auriol2023robustification}. This will be done when designing the output-feedback controller in the next section.}

 Compared to~\citep{deutscher2018output}, our observer does not require assumptions that guarantee the existence of a Byrnes–Isidori normal form for ODE $X_1$. Instead, it relies on the necessary detectability condition stated in Assumption~\ref{Chap_2_ass_detec_PDE}

\section{Output feedback controller} \label{Sec_controller}

In this section, we first design a state-feedback controller to stabilize \eqref{eq_X_0}-\eqref{eq_ODE_distal}. This state-feedback controller will be combined with the state-observer designed in Section~\ref{Sec_observer} to obtain an output-feedback control law. 
Although the target system~\eqref{eq_xi_target}-\eqref{eq_ODE_distal_target_obs} was helpful in designing the state-observer, it may be convenient to transform the coupling term $G_1(x)\alpha(t,1)$ by a term that depends on $\alpha(t,0)$ to design a stabilizing controller. Therefore, before designing the control input, we will first modify the target system~\eqref{eq_xi_target}-\eqref{eq_ODE_distal_target_obs} using a new backstepping transformation.

\subsection{Backstepping transformations and target system}
Consider the following transformations defined by
\begin{align}
\alpha(t,x)&=\check \alpha(t,x)-\int_x^1 \check L(x,y) \check \alpha(t,y)dy\label{Fred_1} \\
\check \alpha(t,x)&=\bar \alpha(t,x)-\int_0^1 \bar L(x,y) \bar \alpha(t,y)dy\label{Fred_2} ,
\end{align}
The kernel~$ \check L$ is a lower triangular matrix  (i.e.,~$( \check L)_{ij}=0$ if~$i< j$) whose components are bounded piecewise continuous functions. The kernel~$\bar L$  is a strictly upper-triangular matrix (i.e.,~$\bar L_{ij}=0$ if~$i \geq j$) whose components are bounded piecewise continuous functions. The transformation~\eqref{Fred_1} is invertible since it is a Volterra transform.  The kernel~$\check L$ satisfies the following set of equations if~$i \geq j$ 
 \begin{align}
&\Lambda^+ \partial_x \check  L(x,y)+\partial_y\check  L(x,y)\Lambda^+=0, \label{eq_L_check_1}\\
&\Lambda^+ \check L(x,x)- \check  L(x,x)\Lambda^+=0, \label{eq_boun_L_check_x}\\
&(\check  L(x,1))_{ij}= (G_1(x)(\Lambda^+)^{-1})_{ij}\nonumber \\
&\quad \quad \quad \quad  +\int_{x}^1 \sum_{k=1}^n \check L_{ik}(x,y) \check G_{kj}(y)\frac{1}{\lambda_j}dy , \label{eq_L_Check_boud}
\end{align}
where the matrix~$\check G(x)$ is strictly upper-triangular (i.e. we have~$\check G_{i,j}(x)=0$ if~$i \geq j$) and satisfies for all~$x \in [0,1]$
\begin{align}
   &(\check  G(x))_{ij}= (G_1(x))_{ij}\nonumber \\
   &+\int_{x}^1 \sum_{k=1}^n \check L_{ik}(x,y) \check G_{kj}(y)dy \quad \text{if~$i < j$}. \label{eq_G_check}
\end{align}

\begin{lem}
 The set of equations~\eqref{eq_L_check_1}-\eqref{eq_G_check}  has a unique solution in~$\mathcal{T}_u$, which is bounded and piecewise continuous. 
 \end{lem}
 \begin{pf}
The proof is a consequence of the triangular structure of the different matrices appearing in the equations. For~$j=1$, equation \eqref{eq_L_Check_boud} can be rewritten as
   ~$$(\check  L(x,1))_{i1}=(G_1(x)(\Lambda^+)^{-1})_{i,1}.$$
Combining this boundary condition with equation~\eqref{eq_boun_L_check_x}, we can solve equation~\eqref{eq_L_check_1} to compute~$\check L_{i1}$ on its domain of definition. For~$j=2$, equation \eqref{eq_G_check} can be rewritten as
$$
(\check  G(x))_{12}= (G_1(x))_{12}+\int_{x}^1  \check L_{11}(x,y) \check G_{12}(y) dy,
$$
which is a Volterra equation that can be solved to obtain~$\check G_{12}(y)$~\citep{yoshida1960lectures}. This in turns gives the kernels~$\check L_{i2}$ using \eqref{eq_L_Check_boud}. Iterating the process allows us to compute the kernel matrix~$\check L$ and the function~$\check G$. 
 \end{pf}

The kernel~$\bar L$ satisfies the following set of equations
\begin{align}
&\Lambda^+ \partial_x \bar  L(x,y)+ \partial_y \bar  L(x,y)\Lambda^+=0,~\bar L(1,y)=0,  \label{eq_L_bar_1}\\
&(\bar L(x,1))_{ij}=(\check G(x)(\Lambda^+)^{-1})_{ij}~ \text{if~$i < j$}, \label{eq_L_bar_2}
\end{align}
Due to its triangular structure, we can obtain a direct expression of~$\bar L$ (and consequently show its existence) using the method of characteristics. The invertibility of the transformation~\eqref{Fred_2} is a consequence of the triangular structure of the kernel~$\bar L$.
The transformations \eqref{Fred_1}-\eqref{Fred_2} map the system  \eqref{eq_xi_target}-\eqref{eq_ODE_distal_target_obs} to the target system
\begin{align}
&\dot{\xi} =\bar A_0\xi+G_3 \bar \alpha(t,1)+G_4X_1(t),\label{Chap_1_eq_X_target_init}\\
&\bar \alpha(t,0)=C_0\xi(t)+Q \beta(t,0)+(Q\gamma_\beta(0)-\gamma_\alpha(0))X_1(t)+  \nonumber\\
&\int_0^1 \bar F^\alpha(y)\bar \alpha(t,y)dy+\int_0^1 F^\beta(y)\beta(t,y)dy+U(t),\\
&\partial_t\bar \alpha(t,x)+\Lambda^+ \partial_x\bar\alpha(t,x) = G_5(x)\bar \alpha(t,0) , \\
&\partial_t \beta(t,x)-\Lambda^- \partial_x\beta(t,x) = G_2(x) \bar \alpha(t,1),\\
&\beta(t,1)=R\bar \alpha(t,1),\quad \dot{X}_1= A_1X_1+E_1\bar \alpha(t,1), 
\label{Chap_1_eq_X_target_bis}
\end{align}
where 
\begin{align*}
\bar F^\alpha(y)=&F^\alpha(y)+\bar L(0,y)+\check L(0,y)+\int_0^1 \check L(0,\nu) \bar L(\nu,y)d\nu \\
&-\int_0^1 F^\alpha(\nu) \bar L(\nu,y)d\nu-\int_0^yF^\alpha(\nu) \check L(\nu,y)d\nu\nonumber \\
&+\int_0^1\int_1^\eta F^\alpha(\nu)\check L(\nu, \eta)\bar L(\eta,y) d\nu d\eta. 
\end{align*}
The upper-triangular matrix function~$G_5$ is defined by 
$$G_5(x)=\bar L(x,0)\Lambda^+ +\int_0^1 \bar L(x,y)G_5(y)dy.$$
We denote $\mathcal{T}_1$ the combination of the transformations \eqref{Fred_1}-\eqref{Fred_2} such that $(\xi, \bar \alpha,\beta, X_1)=\mathcal{T}_1(X_0, u,v,X_1)$. We have the following lemma
\begin{lem}\label{lem_conv_cascade_bis}
If $\xi(t)$, $\bar \alpha(t,0)$ and $X_1(t)$ exponentially converge to zero, then the state $( \xi, \alpha,  \beta,  X_1)$ converges to zero in the sense of the $\chi$-norm. This implies the stabilization of~\eqref{eq_X_0}-\eqref{eq_ODE_distal} .
\end{lem}
\begin{pf}
    The proof is analogous to the proof of Lemma~\ref{lem_conv_cascade}.
\end{pf}

\subsection{Time-delay formulation and state-feedback control law}
Similarly to what has been done to obtain equation~\eqref{Chap_2_eq_delay_alpha_1}, we can adjust the methodology from~\citep{auriol2019explicit} to rewrite $(\bar \alpha(t,0), \xi, X_1)$ as the solution of the following integral differential systems

\begin{align}
    \bar \alpha(t,0)&=\sum_{i=1}^{N} P^\alpha_i  \bar \alpha(t-\tau_i,0) +\int_0^\tau Q_\alpha(\nu)  \bar \alpha(t-\nu,0) d\nu \nonumber \\
    &+C_0\xi(t)+(Q\gamma_\beta(0)-\gamma_\alpha(0))X_1(t)+U(t),\label{eq_neutral_z} \\
    \dot{\xi}(t)&=A_0 \xi+\sum_{i=1}^{n} P_i^\xi  \bar \alpha(t-\tau_i,0)\nonumber \\
    &+ \int_0^\tau Q_\xi(\nu) \bar \alpha(t-\nu,0) d\nu +G_4X_1(t), \label{eq_neutral_X0} \\
    \dot{X}_1(t)&= A_1X_1(t)+\sum_{i=1}^{n} P_i^X\bar \alpha(t-\frac{1}{\lambda_i},0)\nonumber \\
    &+\int_0^\tau Q_X(\nu) \bar \alpha(t-\nu,0) d\nu, \label{eq_neutral_X1}  
\end{align}
where the matrices~$P^\alpha, P^\xi, P^X, Q_\alpha, Q_\xi$ and $Q_X$ can be explicitly computed following the methodology of \citep{auriol2019explicit}. To design a state-feedback control law, we first simplify equation~\eqref{eq_neutral_z} by choosing 
\begin{align}
    U(t)&=\bar U(t)-\sum_{i=1}^{N} P^\alpha_i  \bar \alpha(t-\tau_i,0) -\int_0^\tau Q_\alpha(\nu)  \bar \alpha(t-\nu,0) d\nu \nonumber \\
    &-C_0\xi(t)-(Q\gamma_\beta(0)-\gamma_\alpha(0))X_1(t), \label{Chap_2_eq_U_bar}
\end{align}
so that we have $\bar \alpha(t,0)=\bar U(t)$.
Note that, due to the cancellation of the reflection term~$\sum_{i=1}^{N} P^\alpha_i \bar \alpha(t-\tau_i,0)$, the control law~$U(t)$ will not be strictly proper. Consequently, it may be necessary to low-pass filter it to guarantee the robustness of the closed-loop system~\citep{auriol2023robustification}. Inspired by~\citep{bekiaris2010stabilization}, we define the following change of coordinates
\begin{align}
     \bar X_1(t)=X_1(t)+&\tau^2\int_0^1\big(\int_0^x \mathrm{e}^{- A_1\tau(x-s)} Q_X(\tau(1-s))ds\big)\nonumber \\
     &\bar U(t-(1-x)\tau)dx, \label{Chap_2_transfo_Y_bar} \\
      \bar \xi(t)=\xi(t)+&\tau^2\int_0^1\big(\int_0^x \mathrm{e}^{-A_0\tau(x-s)} \bar Q_\xi(\tau(1-s))ds\big)\nonumber \\
      &\bar U(t-(1-x)\tau)dx, \label{Chap_2_transfo_X_bar}
\end{align}
where $\bar{Q}_\xi(\nu)=Q_\xi(\nu)-G_4\tau \int_0^{\frac{1-\nu}{\tau}}\mathrm{e}^{- A_1(\tau(1-s)-\nu)}$ $Q_X(\tau(1-s))ds$.
We obtain 
\begin{align*}
\dot{\bar{\xi}}(t)&= A_0\bar \xi(t)+E_0\gamma_\beta(0)\bar X_1(t)+\sum_{i=1}^n P_i^\xi \bar U(t-\tau_i+\bar E_0 \bar U(t)\\
     \dot{\bar{X}}_1(t)&=  A_1 \bar X_1(t)+\sum_{i=1}^{n} P_i^X\bar U(t-\frac{1}{\lambda_i})+\bar E_1 \bar U(t)
\end{align*}
where~$\bar E_1=\tau \int_0^1 \mathrm{e}^{-A_1\tau(1-s)}Q_X(\tau(1-s))ds$ and~$\bar E_0=\tau \int_0^1 \mathrm{e}^{-A_0\tau(1-s)}\bar Q_\xi(\tau(1-s))ds$. Let us denote
$$
X_c(t)=\begin{pmatrix}\bar \xi(t) \\ \bar X_1(t)\end{pmatrix},~ A_c=\begin{pmatrix}A_0 &  E_0 \gamma_\beta(0) \\0 &  A_1 \end{pmatrix}.
$$
Inspired by~\citep{artstein1982linear}, we finally define the state~$Z_c$  as
\begin{align}
    Z_c(t)=&X_c(t)+\sum_{i=1}^n \int_{t-\tau_i}^t \mathrm{e}^{A_c(t-s-\tau_i)}\begin{pmatrix}P_i^\xi  \\ 0\end{pmatrix}\bar U(s) ds \nonumber \\
    &+\sum_{i=1}^n \int_{t-\frac{1}{\lambda_i}}^t \mathrm{e}^{A_c(t-s-\frac{1}{\lambda}_i)}\begin{pmatrix}0\\P_i^X \end{pmatrix}\bar U(s) ds .\label{Chap_2_transfo_Z}
\end{align}
Consequently, we obtain~$\dot Z_c(t)=  A_c Z_c(t)+\bar B \bar U(t)$,
with
$$
\bar B=\begin{pmatrix}\bar E_0 \\ \bar E_1 \end{pmatrix}+\sum_{i=1}^n\begin{pmatrix}\mathrm{e}^{-A_c \tau_i}P_i^\xi \\ \mathrm{e}^{-A_c \frac{1}{\lambda_i}}P_i^X \end{pmatrix}.
$$
We are led to the following stabilizability assumption
\begin{assum} \label{Chap_2_ass_stab_PDE}
    The pair~$(A_c,\bar B)$ is stabilizable, i.e., there exists~$K_c$ such that~$A_c+\bar B K_c$ is Hurwitz.
\end{assum}

\begin{thm} \label{th_state_feedback_PDE}
Consider that Assumptions~\ref{Assum_OL_stab} and~\ref{Chap_2_ass_stab_PDE} are verified. 
Consider the control law~$U(t)$ defined by eq.~\eqref{Chap_2_eq_U_bar} where~$\bar U(t)=K_cZ_c(t),$
where~$Z_c$ is defined from~$\bar X$ and~$Y$ using transformations~\eqref{Chap_2_transfo_Y_bar}-\eqref{Chap_2_transfo_X_bar}  and~\eqref{Chap_2_transfo_Z}. Then the closed-loop system~\eqref{eq_X_0}-\eqref{eq_ODE_distal} is exponentially stable.
\end{thm}

\begin{pf}
The exponential stability of~$Z_c$ implies that~$\bar U$ converges to zero. This implies the exponential stability of the state~$\xi$ (using \eqref{Chap_2_transfo_X_bar}) and~$X_1$ (using \eqref{Chap_2_transfo_Y_bar}). This, in turn, implies the exponential convergence of the state~$\bar \alpha(t,0)$ and consequently of~$(X_0, u, v, X_1)$. 
\end{pf}

\subsection{Output-feedback control law}

\begin{thm} \label{Chap_2_th_state_feedback_output_bis}
Consider system~\eqref{eq_X_0}-\eqref{eq_ODE_distal} and that Assumptions~\ref{Assum_OL_stab}, \ref{Chap_2_ass_stab_PDE} and \ref{Chap_2_ass_detec_PDE} are satisfied. Consider the operator  $\mathcal{O}_0$ defined by equation~\eqref{Chap_2_operator_O_0}. Define the observer states $(\hat{X}_0, \hat u,  \hat v, \hat{X}_1)=\mathcal{T}( \hat{\xi},$ $ \hat \alpha, \hat \beta,\hat \nu)$, where $(\hat{\xi}, \hat \alpha, \hat \beta,\hat \nu)$ is the solution of the system~\eqref{Chap_2_eq_xi_target_estim_bis}-\eqref{Chap_2_eq_xi_target_estim_bis_end}. Define the state $(\xi, \bar \alpha, \beta,X_1)=\mathcal{T}_1(\hat X_0,\hat u,\hat v,\hat X_1)$. Define $Z_C$ using transformations~\eqref{Chap_2_transfo_X_bar}-\eqref{Chap_2_transfo_Y_bar}  and~\eqref{Chap_2_transfo_Z}.
 There exists a  low pass filter~$w(s)$ such that the control law defined in the Laplace domain by \begin{align}
    &\hat U(s)=(-\sum_{i=1}^NP_i^\alpha\bar \alpha(t-\tau_i,0)-\int_0^\tau Q_\alpha(\nu)\bar \alpha((t-\nu,0)d\nu  \nonumber \\
    &-C_0\xi(t)-(Q\gamma_\beta(0)-\gamma_\alpha(0))X_1(t)+K_cZ_C(t))w(s) \label{Chap_2_Control_output_final_bis}
\end{align} is strictly proper and exponentially stabilizes~\eqref{eq_X_0}-\eqref{eq_ODE_distal} . 
\end{thm}
\begin{pf}
The closed-loop stability is implied by Theorems~\ref{Chap_2_th_observer_PDE} and~\eqref{th_state_feedback_PDE}. The existence of a low-pass filter that makes the control law strictly proper is shown in~\citep{auriol2023robustification}.
\end{pf}
As explained in~\citep{auriol2023robustification}, having a strictly proper controller guarantees the existence of robustness margins for a broad class of perturbations: input delays, uncertainties on the ODE parameters, uncertainties on the transport velocities~\citep{curtain2012introduction}.

\section{Concluding remarks} \label{Sec_CCl}
In this paper, we have designed an output-feedback controller for a class of \modif{ $n + m$} linear hyperbolic ODE-PDE-ODE systems for which the PDE subsystem is actuated and measured. The proposed approach combines backstepping transformations and a rewriting of the original system as a time-delay system. We first designed a state observer that was then combined with a state-feedback control law to obtain an output-feedback controller. In future works, we will consider networks of ODE and PDEs with a more complex structure (star-shaped networks).

\bibliography{Biblio}

\end{document}